\documentclass{article}

\usepackage{amssymb}
\usepackage{amsmath}

\usepackage[utf8]{inputenc}
 \newmuskip\pFqmuskip

\newcommand*\pFq[6][8]{%
  \begingroup % only local assignments
  \pFqmuskip=#1mu\relax
  \mathcode`\,=\string"8000
  \begingroup\lccode`\~=`\,
  \lowercase{\endgroup\let~}\pFqcomma
  {}_{#2}F_{#3}{\left[\genfrac..{0pt}{}{#4}{#5};#6\right]}%
  \endgroup
}
\newcommand{\pFqcomma}{\mskip\pFqmuskip}

\title{Applications of Laguerre transform to solve Schr\"{o}dinger-type equations and Differential Equations of order  four}
\author{Gabriel L\'{o}pez Garza\\ Universidad Aut\'{o}noma Metropolitana Iztapalapa\\
Ciudad de México, México}

\begin{document}
\maketitle

\section*{Abstract}
The finite Laguerre transform is applied to solve Differential Equations Problems of order higher than two and a one-dimensional steady-state Schr\"{o}dinger equation, by using elementary Linear Algebra methods.

%% %%%%%%%%%%%%%%%%%%%%%%      INTRODUCTION
\section{Introduction}

The theory of finite Sturm-Liouville transform, as studied in  \cite{CE}, has been applied to solve differential and partial differential equations since the fifties of the last century.
McCully, in particular, developed the theory of  Laguerre transform in \cite{McC}, where the author after calculating some transforms formulas, employs the finite Laguerre transform for solving the one-dimensional nonhomogeneous heat equation as an instance of an application.
Since then, the Laguerre transform has become a useful tool for solving differential equations and for other applications.
Taking into account the large body of literature related to Laguerre transforms it is fair to suggest that wherever Laguerre polynomials appear in relation to ordinary and partial differential equations, it is possible to use the finite Laguerre transform.
The present article could be inscribed in support of such a hypothesis.
In fact, related to this, it can be mentioned the case of the article \cite{Al} where  recently, Alhaidari utilizes Laguerre polynomials to solve, among others, some examples of the Sch\"{o}dinger equation. 
To this aim, he applies the basic properties of Laguerre polynomials, as well as the theory of three-terms recursion sequences for  difference equations of order two, as presented in \cite{Ko}.
Approximations to the solutions of the Schr\"{o}dinger equation for the Morse potential case, in terms of the Laguerre polynomials, are known \cite{G}, \cite{V}, but the technique introduced in \cite {Al} seems to be new.
Nevertheless, Alhaidari never applies explicitly the Laguerre finite transform which, as the present article shows, simplifies and clarifies the pertinence of such kinds of techniques. 
In particular, the relevance of the most important Sturm-Liouville transform $\mathcal{T}[L[y_n(x)]]=\lambda_n \mathcal{T}[y_n(x)]$, where $L$ is a second order differential operator and $\lambda_n,y_n(x)$ are the eigenvalues and eigenfunctions related to the operator respectively, is completely hidden in Alhaidari's work.
 By adding and subtracting appropriate terms, the author avoids the use of transforms, which is valid, but such procedures do not permit the implementation of algorithms efficiently.
By employing the finite Laguerre transform method, such difficulties can be overcome without great difficulty.

Another application treated in this article deals with the iteration of a Strum-Liouville operator, which transform is  
{\small $\mathcal{T}[L^2[y_n(x)]]=(\lambda_n)^2\mathcal{T} [y_n(x)],$ } simply,
but when used in physical applications it gives rise to equations that look very complicated and difficult to solve, but with finite transform methods it is not so.
In such class of problems may be included the Laguerre type orthogonal polynomials discovered by Krall  \cite{Kr}, \cite{Kr2}, which are instances of eigen-solutions of  fourth-order differential equations that do not satisfy second order differential equations.

The finite transforms method associates with the equations studied in this article systems of linear equations which in the case of the Schr\"{o}dinger equation leads can be reduced to a three-term recursion formula which are completely classified and known \cite{Ko}. 
In the case of equations of order bigger than two it is not possible to reduce the equations to systems of three-terms, but it is possible to reduce them two a triangular system which somehow, it is easy to solve.
In both last cases, ordinary differential equations solutions are solved by elementary linear algebra methods.

In section \ref{mathset}, the mathematical frame of the Sturm-Liouville finite transform is briefly summarized to set the context in which the Laguerre finite transform is inscribed. 
The basic properties of the Laguerre transform are also exposed in this section.
 In section \ref{EAPP} are solved the Coulomb and Morse instances of the Schr\"{o}dinger equation as well as a Laguerre-type equation of fourth-order.
The literature on Laguerre transforms is immense and, hence, it is difficult to know if some transform is already known. 
The criteria followed in this article is that, if a formula is not in the classical book of  Erd\'{e}lyi \cite{Er} or in  \cite{Be}, or if a Laguerre transform is not in \cite{P}, the proof must be exhibited.
In the Appendix section \ref{app}, the reader can find the proofs of some of such Laguerre transforms. 
That the solution of the fourth-order equation studied in this article is a linear combination of Laguerre polynomials is also proved in the Appendix.

%%%%%%%%%%%%%%%%%%%%%%%%%%%    SETTING
\section{Mathematical setting}
\label{mathset}

%%%%%%%%%%%%%%%%%%%%%%    DISCRETE TRANSFORMS
\subsection{Sturm-Liouville Finite Transform}
\label{SLft}
Consider a  Suturm-Liouville boundary value problem
\begin{eqnarray}
\label{slp} \mathcal{L}[y](x)=\frac{d}{dx}\left(p(x)\frac{d}{dx}y(x)\right)-q(x)y(y)=-\lambda r(x) y(x)\\
\label{slbv} a_1y(\alpha)+b_1y'(\alpha)=0,\qquad a_2y(\beta)+b_2y'(\beta)=0,
\end{eqnarray}
when the eigenfunctions are polynomials.
 Specifically, in this article is consider the case of Laguerre polynomials 
which, as the reader may recall, are solutions of
\begin{eqnarray}
\label{Laguerre} \mathcal{L}[y](x)&=&\frac{d}{dx}\left(x^{\nu+1}e^{-x}\frac{d}{dx}y(x)\right)=-n x^\nu e^{-x}y(x),0<x<\infty.
\end{eqnarray}
And given that $\beta=+\infty$, the boundary condition ~(\ref{slbv}) become
\begin{eqnarray}
\label{modbv} 
 \lim_{\beta\to+\infty}  r(\beta)(a_2y(\beta)+b_2y'(\beta))=0,
\end{eqnarray}
where $r(x)= x^\nu e^{-x}$, accordingly.
Solutions of (\ref{Laguerre}) are denoted $L^\nu_n(x)$, for $\nu\neq 0$ and, as is usually done in the literature, $L^0_n(x)$ is denoted $L_n(x)$ simply. The number $\nu$ is called the order of the corresponding Laguerre polynomials.

Solutions of the equation (\ref{slp}) which satisfy boundary conditions  (\ref{slbv}), denoted by $P_n(x), \,n\geq 0$. It is well known that suitable differential functions $f(x)$ (i.e. functions satisfying the respective boundary value conditions and sufficiently differentiable to satisfy the corresponding differential equation), have an expansion, called Sturm -Liouville expansion, of the form
\begin{eqnarray}
\label{slt} f(x) &= & c_0 P_0(x)+c_1P_1(x)+\cdots=\sum_{m=0}^\infty c_mP_m(x),\\
\nonumber \mbox{where } c_m &=&\frac{\int_\alpha^\beta r(x)f(x)P_m(x)dx}{\int_\alpha^\beta r(x)P_m^2(x)dx}.
\end{eqnarray}
Of course, assuming $r(x)>0,\alpha\leq x \leq \beta$,  to be continuous in the interval of definition, the integral $\int_\alpha^\beta r(x) h(x)g(x)dx$ defines an inner product $\langle h,g\rangle_r$ in a corresponding space of integrable functions.

From formula~(\ref{slt}) it is possible to establish a correspondence $f(x)\leftrightarrow \{c_m\},m\geq 0$, so that
the sequence $\{c_n\}$  is known as the Sturm-Liouville finite transform of $f(x)$, which is usually denoted by
\begin{eqnarray}
\label{transformf} \mathcal{T}[f](s)=\{c_m\}\stackrel{def}{=} c_0+c_1s+c_2s^2+\cdots
\end{eqnarray}
With the purely formal notation defined in~(\ref{transformf}) we mean the sum and Cauchy product of sequences, where $s^n$ is the sequence
\begin{eqnarray}
s^n\stackrel{def}{=} \{\underbrace{0,0,\dots,0}_{n},1,0,0,\dots\}
\end{eqnarray}
 The sum of sequences $\{a_n\},\{b_n\}$ is defined as usual by $\{a_n\}+\{b_n\}\stackrel{def}{=}\{a_n+b_n\}$ and the Cauchy product of sequences  is defined by the sequence
\begin{eqnarray}
\label{CP} \{a_n\}\{b_n\}=\{a_0b_0, a_0b_1+a_1b_0,a_0b_2+a_1b_1+a_2b_0,\dots,\sum_{\tau=0}^na_\tau b_{n-\tau},\dots\}.
\end{eqnarray}
 The sequence $\{k,0,0,\dots\}$ is simply written
$$k\stackrel{def}{=}\{k,0,0,\dots\}=k\{1,0,0,\dots\},$$
 as a simplification in the notation  currently used,
so by (\ref{transformf}) it is understood
\begin{eqnarray*}
  \mathcal{T}[f](s)=\{c_m\}&= &c_0+c_1s+c_2s^2+\cdots+c_ns^n+\cdots\\
 &=& \{c_0,0,\dots\}\{1,0,\dots\}+ \{c_1,0,\dots\}\{0,1,0,\dots\}+\cdots \\
 & &+\{c_n,0,\dots\} \{\underbrace{0,0,\dots,0}_{n},1,0,0,\dots\}+\cdots.
\end{eqnarray*}
Observe that convergence is not germane to this notation, and (\ref{transformf}) merely states that the term $c_m$  of the sequence occupies the same place as the $1$ in the sequence $s^m$.
Finally, the last step in solving differential equations via transforms is to find the inverse transform. In the case of finite transforms, the inverse transform is find simply writing down the corresponding series in terms of eigenfunctions \cite{CE}, i. e. series (\ref{slt}) gives directly the inverse transform of a sequence 

\begin{eqnarray}
\label{slti}  \mathcal{T}^{-1}[\{c_n\}]=c_0P_0(x)+c_1P_1(x)+c_2P_2(x)+\cdots
\end{eqnarray}
Of course, convergence matters in this case, but this aspect is well-known for Laguerre polynomials or at least boundedness of solutions is warranted for the problems studied in this article (see \cite{Al} for boundedness of solutions of the Sch\"{o}dinger equation).
Nevertheless, such expressions depend on the characteristics of each particular examples which have to be specified. %To illustrate this particularities of the theory some examples will be presented in the following subsections.

%%%%%%%%%%%%%%%%%%%%%%%%   LAGUERRE TRANSFORMS of orden \nu
\subsection{ Laguerre transforms of order $\nu$}
The Laguerre polynomials are defined by 
\begin{eqnarray*}
L_n^\nu (x)= \sum_{m=0}^n \binom{n+\nu}{n-m} \frac{(-x)^m}{m!},
\end{eqnarray*}
where $\displaystyle{\binom{\alpha}{\beta}=\frac{\Gamma (\alpha +1)}{\Gamma (\beta +1)\Gamma(\alpha-\beta+1)}.}$
The basic orthogonality relation between the Laguerre polynomials is given by 
\begin{eqnarray}
\label{ort} \int_0^\infty e^{-x}x^\nu L_n^\nu(x) L_m^\nu (x)dx=\frac{\Gamma(n+\nu+1)}{n!}\delta_{nm},
\end{eqnarray}
where $\delta_{nm}$ is the Kronecker delta symbol.
Such relation leads to the definition of the Laguerre transform of order $\nu$:
\begin{eqnarray}
\label{deftrans} \mathcal{T}^\nu[f(x)]=\left\{\int_0^\infty e^{-x}x^\nu L_n^\nu(x)f(x)dx\right\}=\{c^\nu_n\},
\end{eqnarray}
it must be emphasize that the Laguerre transform is a sequence of numbers in $\mathbb{C}$.
%In section \ref{scs} Laguerre polynomials are  normalized by multiplying by  $\sqrt{\frac{n!}{\Gamma(n+\nu+1)}}$.
The inverse Laguerre transform is defined by
\begin{eqnarray}
\label{invtll} f(x)= \mathcal{(T^\nu)}^{-1}[\{c_n^\nu\}]=\sum_{k=0}^\infty c^\nu_nL_n^\nu(x).
\end{eqnarray}

%%%%%%%%%%%%%%%%%%%%%%%%%      EXAMPLES
\section{Examples of applications}
\label{EAPP}

%%%%%%%%%%%%%%%%%%%%%%%%%%%%  SCHRODINGER

\subsection{Applications to the Schr\"{o}dinger equation}
\label{scs}

In this section is consider the equation
\begin{eqnarray}
\label{sch1}- \frac{1}{2}\frac{d^2 }{dr^2}\psi (r)-(V(r)+E)\psi (r)=0
\end{eqnarray}
which in appropriate units ($\hbar$=M=1) is the steady state Sch\"{o}dinger equation defined in a one dimensional space, where $V(r)$ is a potential function and $E$ is the energy. Under the change of coordinates (see \cite{Al}), given by $\lambda^{-1}\xi(x)=d x/dr$, where $\lambda\geq 0$ has inverse length units,   equation (\ref{sch1}) becomes
\begin{eqnarray}
\label{54} \lambda^2\xi^2 \left[ \frac{d^2 }{dx^2}\psi(x)+\frac{1}{\xi}\frac{d \xi}{dx}\frac{d}{dx}\psi(x)-\frac{2}{\lambda^2\xi^2 }W(x)\psi(x)\right]=0,
\end{eqnarray}
where $W(x)=V(r)-E$.
To obtain  a Laguerre-type equation the change of coordinates must satisfy $x(r)\geq 0$ and setting $\frac{1}{\xi} \frac{d \xi}{dx}=\frac{a}{x}$, leads to $\xi(x)=x^ae^{bx}$.

% ñ hoy
In this way, equation (\ref{54}) becomes

\begin{eqnarray}
\label{sch2p}\lambda^2 \xi^2 \left[\frac{d^2}{dx^2}\psi(x)+\left(\frac{a}{x}+b\right)\frac{d}{dx}\psi(x)+\left(A_{+}+\frac{A_{-}}{x^2}-\frac{A_0}{x}\right)\psi(x)\right]=0,
\end{eqnarray}
where $A_{\pm}, A_0,a, b$ are real parameters determined in terms of $V(r)$ and  $E$.

To solve equation (\ref{sch2p}) it is proposed a solution of the form
\begin{eqnarray}
\label{solsch} \psi(x)= x^{\alpha}e^{-\beta x}y(x),
\end{eqnarray}
where $y=\sum_{k=0}^\infty c_kL_k^\nu(x)$, and $L_n^\nu(x)$ are the Laguerre polynomials of order $\nu$, and $\alpha,\beta, \nu$ are dimensionless parameters, free for the moment, but to be determined according to the concrete examples to be solved below.

To solve (\ref{sch2p}) the use of the finite Laguerre transform is introduced.
Many of the following transforms are known \cite{P} or are obtained by direct calculation by using Laguerre polynomial properties found  in \cite{Be} or in \cite{Er}.

The following formulas (\ref{tnu3}), (\ref{tnu4}) and (\ref{tnu5}), are proved  in section \ref{app},
since they are not found in any of the items in the bibliography: \cite{Be}, \cite{P} or \cite{Er}.

%%%%%%%%%%%%%%%%%%%%%%%%%%    THEOREM 2    %%%%%%%%%%%%%%%%%%%%%%%%
\noindent {\bf Theorem 1.}
Let $\mathcal{T}^\nu[f(x)]=\{c_n^\nu\}$, then

\begin{eqnarray}
\label{tnu1}\mathcal{T}^\nu\left[\frac{d}{dx}f(x)\right]& = &\left\{ \sum_{k=0}^{n} c_k^\nu -\nu \sum_{k=0}^nc_k^{\nu-1}\right\}\\ %\left\{ c_n^\nu -\nu \sum_{k=0}^nc_k^{\nu-1}+\sum_{k=0}^{n-1} c_k^\nu\right\}, \\
\label{tnu2}\mathcal{T}^\nu\left[\frac{f(x)}{x}\right]   & = &  \left\{ \sum_{k=0}^n c^{\nu-1}_k\right\},\\
 \label{tnu3}\mathcal{T}^\nu\left[   xf(x)          \right]                      & = &  \left\{ (2n+\nu+1)c_n^\nu-(n+1)c_{n+1}^\nu-(n+\nu)c_{n-1}^\nu\right\},\\
\label{tnu4} \mathcal{T}^\nu\left[ x\frac{d}{dx}f(x)\right] & = &\{nc_n^\nu-(n+1)c_{n+1}^\nu\}             , \\ 
\nonumber     \mathcal{T}^\nu\left[ x\frac{d^2}{dx^2}f(x)\right] & = & \left \{ -(\nu+1)\left(\sum_{k=0}^{n}c_k^\nu-\nu\sum_{k=0}^nc_k^{\nu-1}\right)
                                                     -(n+1)c^\nu_{n+1}\right\} .\\
\label{tnu5}
\end{eqnarray}

%The proof of  theorem 2  is in section \ref{app}.
To begin with, the Coulomb problem is solved next.

\subsubsection{Coulomb problem}
%\noindent {\it Example, the Coulomb problem $a=b=0$.}

 Setting the parameters $a=b=0$ follows that $\xi(x)=1$, so that $x=\lambda r$. The corresponding form of equation (\ref{sch2p}) for this problem  is
%In this way, equation (\ref{54}) becomes

\begin{eqnarray}
\label{sch2}\frac{\lambda^2 \xi^2}{x} \left[x\frac{d^2}{dx^2}\psi(x)%+(a+bx)\frac{d}{dx}\psi(x)
+\left(A_{+}x+\frac{A_{-}}{x}-A_0\right)\psi(x)\right]=0,
\end{eqnarray}
where $A_0=\frac{2Z}{\lambda}, A_-=-\mathcal{l}(\mathcal{l}+1), A_+=\frac{2E}{\lambda ^2}$ .

To solve equation (\ref{sch2}) it is proposed a solution of the form
\begin{eqnarray}
\label{solsch} \psi(x)= x^{\alpha}e^{-\beta x}y(x),
\end{eqnarray}
where $y=\sum_{k=0}^\infty c_kL_k^\nu(x)$, $L_n^\nu(x)$ are the Laguerre polynomials of order $\nu$, and $\alpha,\beta, \nu$ are dimensionless parameters, free for the moment, but to be determined according to the concrete examples to be solved below.

%%%%%%%%%%%%%%%%%%      HOY
By substituting $\psi$ of (\ref{solsch}) in (\ref{sch2}) the following equation is obtained, after canceling factors:
\begin{multline}
\label{sch3}  x\frac{d^2}{dx^2}y(x)+(2\alpha -2\beta x)\frac{d}{dx}y+(\beta^2+A_{+})xy(x)+\\
             +[\alpha^2-\alpha+A_-]\frac{y(x)}{x}-(A_0+2\alpha \beta)y(x)=0.
\end{multline}
%  end ñ hoy

It is possible to find a solution of equation (\ref{sch3}) of the form $\psi(x)=x^\alpha e^{-\beta x}y(x)$, where $y(x)=\sum_{k=0}^\infty c_kL_k^\nu(x)$ since $c_n^\nu$ can be calculated  and hence, solve equation (\ref{sch3}), by means of the finite transform method.
%%%%%%%%%%%%%%%%%%%%%%%%    HOY

Applying formulas (\ref{tnu1}) to (\ref{tnu5}) to equation (\ref{sch3}) it is possible to solve the Coulomb problem. The complete  transform of equation (\ref{sch3}) is  %presented in the appendix, section \ref{app}. 
\begin{multline}
\label{buena}
  [-\nu-1+2\alpha-2\beta n+(\beta^2+A_+)(2n+\nu+1)-\\
-(A_0-2\alpha) c_n^\nu   -(\beta^2+A_+) (n+\nu) c_{n-1}^\nu-\\
-[(n+1)(1-2\beta+\beta^2+A_+)]c_{n+1}^\nu+\\
+( \nu+1-2\alpha)\left(-\sum^{n-1}_{k=0}c_k^\nu+\nu\sum_{k=0}^{n}c_k^{\nu-1}\right)+\\
+(\alpha^2-\alpha+A_-)\sum_{k=0}^n c_{k}^{\nu-1}=0.
\end{multline}

\noindent{\bf Remark.} Note that equation  (\ref{buena}) is not a three-terms recursion formula since it does include terms order $\nu-1$ as well as terms of order $\nu$. Nevertheless it is possible to solve some equations by restricting some of the coefficients of formula (\ref{buena}) by considering physical parameters.

The parameters $\alpha$, $\beta$, $\nu$ are chosen in (\ref{buena}) in such a way that the equation does not contain terms of order $c^{\nu -1}_k$ nor terms that include a summation sign.
 In this procedure, it is valid to appeal to %the Criterion of Simplicity
Theorem 6.1 \cite{Ko}[p. 139 case II], since, as already mentioned, $\alpha, \beta, \nu$ are free parameters. Then it is valid to choosing
\begin{eqnarray}
\label{al} 0 & = & \alpha^2-\alpha +A_-\\
\label{be}\beta & = & \frac{1}{2}\\
\label{nu} \nu &=& 2\alpha -1.
\end{eqnarray}
From (\ref{al}) it is obtain  $\displaystyle{\alpha=\frac{1\pm\sqrt{1-4A_-}}{2}}$,
hence by (\ref{nu}) $\nu^2=1^2-4A_-$, so that $\alpha=(1\pm \nu)/2$. Therefore, for Coulomb problem $\alpha=(1\pm\sqrt{1-4A_-})/2$, $\nu=\pm \sqrt{1-4A_-}$, and
$1/4\geq A_-$. %  After taking transforms equation (\ref{sch2}) becomes (\ref{buena})  (see appendix section \ref{app}).
Finally, the choice of parameters
\begin{eqnarray}
\label{param} A_0=\frac{2Z}{\lambda},\qquad A_-=-\mathcal{l}(\mathcal{l}+1),\qquad A_+=\frac{2E}{\lambda ^2},
\end{eqnarray}
 leads to 
\begin{eqnarray}
\label{paramc}\alpha= \mathcal{l}+1, \qquad \nu = 2\mathcal{l}+1,\qquad % \beta=\frac{1}{2}
\end{eqnarray}
where $Z$ is the electric charge and $\mathcal{l}$ is the angular momentum quantum number. % and, taking into consideration physically acceptable constants (see \cite{Al}[sec. A1]),

Equation  (\ref{buena}) leads to a three terms recurrence relation as follows. 
 By substituting (\ref{al}), (\ref{be}), and (\ref{nu}) in the transformed equation (\ref{buena}) it is obtain after simplification

\begin{multline}
\label{CouT}    (n+1)c_{n+1}^\nu+\left(\frac{n+A_0 +\alpha}{ \left(A_++\frac{1}{4}\right)}-(2n+\nu+1)\right)c_n^\nu
+(n+\nu)c_{n-1}^\nu=0.
\end{multline}
Next, substitute the values of (\ref{param}) in (\ref{CouT}). After collecting terms and simplification it follows
\begin{eqnarray}
\label{CouT2} (n+1)c_{n+1}^\nu-2\left(  \frac{-Z}{4E}\sin \phi + (\left( n+\alpha) \cos\phi \right)  \right)c_n^\nu+(n+2\alpha-1)c_{n-1}^\nu=0,
\end{eqnarray}
where $\displaystyle{\cos \phi =\frac{4(2E)-\lambda^2}{4(2E)+\lambda^2}}$.
As noticed by  Alhaidari \cite{Al}, equation  (\ref{CouT2}) is  related to the Meixner-Pollaczek  three-terms recurrence relation in \cite{Ko}[eq. (9.7.3), p. 213].
 Note that (\ref{CouT2}) is not identical with the Meixner-Pollaczek formula since the order $\nu=2\alpha -1$ differs from the order $(\lambda)$  in equation (9.7.3) in \cite{Ko}. Consequently, by setting $\lambda=\frac{\nu+1}{2}=\mathcal{l}+1$  a solution of equation (\ref{sch2}) is found by choosing the coefficients $c_n^\nu$ as % &
\begin{eqnarray}
\label{MP}    c_n^\nu=P_n^{\left(\frac{\nu+1}{2}\right)}(z;\phi)=\frac{\left(\frac{\nu+1}{2} \right)_n}{n!}e^{in\phi} \pFq{2}{1}{-n,\lambda+iz}{2\lambda}{1-e^{-2i\phi}},
\end{eqnarray}
where, $_2F_1[\cdot]$ is a hipergeometric function, $z=-Z/2E$, and $(k)_n$ is the Pochhammer symbol (for definitions and notation see \cite{Ko}[Ch. 1 sec. 1.4], for instance).
Therefore, taking into account the relations in (\ref{paramc}) and substituting in  (\ref{MP}), the solution of equation (\ref{sch2})  is
\begin{eqnarray}
\label{solcoul}\psi (x)=x^{\mathcal{l}+1}e^{-x/2 }\sum_{k=0}^\infty P_k^{\mathcal{l}+1}(z;\phi)\,L_k^{2 \mathcal{l}+1}(x),
\end{eqnarray}
 which is equivalent to equation (60) in \cite{Al}.
%where $\nu=1+\sqrt{1+4\mathcal{l}(\mathcal{l}-1)}$. 

\noindent{\bf Remark.} Recall that  in this paper the Laguerre polynomials $L_n^\nu(x)$ are not normalized and that $z,\phi$ are fixed parameters in  (\ref{solcoul}).

%%%%%%%%%%%%%%%%%      Morse oscillator

\subsubsection{One dimensional Morse oscillator}

For this case the values of $a,b $, are $a=1,\;b=0$ so that $\xi(x)=x$. Hence by substituting $\xi,\xi '$ in equation (\ref{sch2p}) and after simplification it is obtain
\begin{eqnarray}
\label{Moc1} x \frac{d^2}{dx^2}\psi(x)+\frac{d}{dx}\psi(x)+\left( A_+x-A_0+\frac{A_-}{x}\right)\psi(x) =0.
\end{eqnarray}
Now, again,  it is proposed a solution $\psi(x)=x^\alpha e^{-\beta x} y(x)$ of (\ref{Moc1}), which after substitution gives
\begin{multline}
\label{schMo}  x\frac{d^2}{dx^2}y(x)+(1+2\alpha-2\beta x)\frac{d}{dx}y(x)+
\\+\left( (\beta^2+A_+)x-(\beta+2\alpha+A_0)+\frac{\alpha^2+A_-}{x} \right)y(x)=0.
\end{multline}
After taking transforms in (\ref{schMo}) with formulas  (\ref{tnu1}) to (\ref{tnu5}), follows the next equation for the Laguerre coefficients of $y(x)=\sum_{k=0}^\infty c_n^\nu L_n^\nu(x)$:

\begin{multline}
\label{solTMSS} c_n^\nu[n(A_+-2\beta+2\beta^2)+\beta(-1-2\alpha+\beta+\beta\nu)+2\alpha-\nu(1-A_+)-A_0+A_+ ]+\\
-c_{n+1}^\nu [(n+1)((1-\beta)^2+A_+)]-\\
-c_{n-1}^\nu(\beta^2+A_+)(\nu+n)+\\
+(-\nu+\alpha(2+\alpha)+A_-)\sum_{k=0}^{n-1}c_k^\nu+\\
+(\nu^2-2\nu\alpha+\alpha^2+A_-)\sum_{k=0}^n c_k^{\nu-1}=0.
\end{multline}

With the same simplification criterion as in the Coulomb example, by appealing to Theorem 6.1 \cite{Ko}[p. 139 case II], it  is obtained the following system by equating the coefficients of the sum symbols to zero
\begin{eqnarray}
\nonumber -\nu+2\alpha & = & 0\\
\label{sysM}\alpha^2+A_-&  = & 0,
\end{eqnarray}
from (\ref{sysM}) it follows that $\nu=2\alpha, \alpha^2=-A_-$, and, moreover, it is possible to set $\beta =1/2$ too, as in the Coulomb example. On the other hand the values of $A_0, A_-, A_+$ given in \cite{Al} are
\begin{eqnarray}
\label{parM} A_0=\frac{ -2V_1}{\lambda ^2},\quad  A_-=\frac{2E}{\lambda^2},\quad  A_+=\frac{-2V_2}{\lambda ^2}.
\end{eqnarray}
By substituting the values of $\alpha,\nu,\beta,A_0,a_-,A_+$ in (\ref{solTMSS}) a three-term recursion system is obtained:
\begin{multline}
\label{solFM} (n+1)c_{n+1}^\nu-c_n^\nu\left[  2n\frac{A_+ -\frac{1}{4}}{A_++\frac{1}{4}}-\frac{\frac{1}{2}(1+\nu)+A_0}{A_+-\frac{1}{4}}+\nu+1\right]+(n+\nu)c_{n-1}^\nu=0.
\end{multline}
Setting $$\cos \phi =\frac{A_+ -\frac{1}{4}}{A_++\frac{1}{4}}$$
equation (\ref{solFM}) becomes
\begin{eqnarray}
\label{finM} (n+1)c_{n+1}^\nu-2c_n^\nu\left[ \left( n+\frac{\nu+1}{2}\right)\cos \phi+\frac{A_0}{2\sqrt{A_+}}\sin \phi\right]+(n+\nu)c_{n-1}^\nu=0.
\end{eqnarray}
Now, formula (\ref{finM}) can be compared with Meixner-Pollaczek recursion formula in \cite{Ko}[eq. (9.7.3), p. 213]:
\begin{multline}
\label{MPM}  (n+1)P_{n+1}^{(\omega)}\left( \frac{A_0}{2\sqrt{A_+}};\phi\right)-\\
-2\left[ \frac{A_0}{2\sqrt{A_+}}\sin \phi+ ( n+\omega)\cos \phi\right] P_n^{(\omega)}\left( \frac{A_0}{2\sqrt{A_+}};\phi\right)\\
+(n+\nu)P_{n-1}^{(\omega)}\left( \frac{A_0}{2\sqrt{A_+}};\phi\right)=0. 
\end{multline}
Set $\omega =(\nu+1)/2$ then the solution of the one dimensional Morse operator is
\begin{multline}
\label{finsolM} \psi(x)= x^{\nu/2}e^{-x/2}\sum_{k=0}^\infty\frac{(\nu+1)_k}{k!}e^{-in\phi} \pFq{2}{1}{-n,,\frac{\nu+1}{2}+i\frac{A_0}{2\sqrt{A_+}}}{\nu+1}{1-e^{-2i\phi}}L_k^\nu(x),
\end{multline}
where $\nu=\frac{2\sqrt{-2E}}{\lambda}$ and $(\nu+1)_k=(\nu+1)(\nu+2)\cdots (\nu+k-1).$ \hfill$\square$

\noindent {\bf Remark.} Given the first argument $-n$ in the hypergeometric  function %$\pFq{2}{1}{\phantom{\cdot},\phantom{\cdot}}{\phantom{\cdot}}{\phantom{\cdot}}$ 
$_2F_1[\cdot]$ in (\ref{finsolM}) it is worth to mention that such function is a polynomial for each $n$, since $(-n)_{n+1}=0$.

%%%%%%%%%%%%%%%%%%    END SCHRODINGER

%%%%%%%%%%%%%%   SECOND EXAMPLE 

\subsection{Second example, Laguerre-type equations of order four}
The problem (\ref{La4}) shown below has been of interest in the literature \cite{Kr}, \cite{Kr2}. The following Theorem 2, shows that it is possible to solve fourth-order Laguerre-type equations with finite transform methods as follows. Note that $L_n^0(x)$ is simply denoted by $L_n(x)$ and similarly $\mathcal{T}^0=\mathcal{T}$.

\noindent{\bf Theorem 2.}{\it  The problem,
\begin{eqnarray}
\label{La4}
\begin{cases}
\displaystyle{\frac{d^2}{dx^2}\left(x^2e^{-x}\frac{d^2 y}{d x^2}\right) -\frac{d}{dx}\left(( [2R+2]x+2)e^{-x}\frac{dy}{dx}\right) =e^{-x}\lambda_m y,}\\
  y(0)  =R>0,\quad y'(0)=\frac{\lambda_m y(0)}{-2R}=-\lambda_m =-m(m+2R+1),\, m\in\mathbb{N},
\end{cases}
\end{eqnarray}
has solutions of the form 
\begin{multline}
\label{solLag} y_m(x)=-L_0(x)-L_2(x)-\cdots-L_{m-1}(x)+(R+m)L_m(x),\quad 0\leq x<\infty
\end{multline}
where 
$$L_n(x)=\sum_{k=0}^n(-1)^k\frac{n!}{(n-k)!(k!)^2}x^k,\;n=0,1,\dots$$
are the order zero Laguerre polynomials of degree $n$.}

The proof of theorem 2 is in section \ref{app}.

\noindent {\bf Remark. } It is worth to mention that solutions obtained by Krall  in \cite{Kr}, \cite{Kr2} for problem (\ref{La4}) are of the form 
\begin{eqnarray}
\label{ksol} R_n(x)=\sum_{k=0}^n\frac{(-1)^k}{(k+1)!}{n\choose k}[k(R+n+1)+R]x^k, \;n=0,1,\dots .
\end{eqnarray}
Such solutions (\ref{ksol}) are obtained by the Frobenius method.
That the solutions obtained by Krall are equal to the solution s obtained in this article i. e., that $y_m(x)=R_m(x),\, m=0,1,\dots$ for all $x\in [0,\infty)$ is shown by induction in the appendix \ref{app}. The novelty of the solutions $y_m(x)$ of (\ref{La4}) is the method of obtaining them that is described below.

\noindent {\bf Method to obtain solutions of (\ref{La4}).} Using (\ref{Laguerre})  the differential equation in (\ref{La4}) may be written as
\begin{eqnarray}
\label{La41} \mathcal{L}^2[y]-(2R+1)\mathcal{L}[y]+2\frac{dy}{dx}-2\frac{d^2y}{dx^2}=\lambda_my,
\end{eqnarray}
where $L$ is defined in identity (\ref{Laguerre}).
The Laguerre transforms of  $ \mathcal{T}[L[y]]$, $T[L^2[y]]$, $ \mathcal{T}[\frac{dy}{dx}]$ are known \cite{McC}, and are included below for the convenience of the reader:
\begin{eqnarray}
\label{tlLa}  \mathcal{T}[\mathcal{L}[y]]& = &\{-nc_n\},\mbox{ where }  \mathcal{T}[y]=\{c_n\},\;n=0,1,2,\dots\\
\label{tlLa2}  \mathcal{T}[\mathcal{L}^2[y]]&=&\{n^2c_n\},\\
\label{tlLad}  \mathcal{T}\left[\frac{dy}{dx} \right]&=&\left\{ \sum_{k=0}^nc_k-y(0)\right\}.%=\sigma (\{c_n\}-\{y(0)\}).
\end{eqnarray}
By integration by parts as in \cite{McC}, or iteration in formula (\ref{tlLad}), it is possible to obtain the formula
\begin{eqnarray}
\label{tlLad2}   \mathcal{T}\left[\frac{d^2y}{dx^2} \right]&=& \left\{ \sum_{k=0}^n (k+1)c_{n-k}-(n+1)y(0)-y'(0)\right\}%\\ 
% \nonumber                  & = &                     \sigma^2(\{c_n\}- y(0))-y'(0)\sigma ,
\end{eqnarray}
where $y'(0)=\frac{dy}{dx}\big|_{x=0}$.
After formulas (\ref{tlLa}), (\ref{tlLa2}), (\ref{tlLad}), and (\ref{tlLad2})  by taking transforms in equation (\ref{La41}), it is  obtain
\begin{multline}
\label{LaT1}                           \{n(n+2R+1)c_n\}+2\left\{ \sum_{k=0}^nc_k-y(0)\right\}-\\
-2 \left\{ \sum_{k=0}^n (k+1)c_{n-k}-(n+1)y(0)-y'(0)\right\} =   m(m+2R+1)\{c_n\}
\end{multline}
%\begin{multline}
%\label{LaT1}2 \sigma\{c_n- y(0)\}-2( \sigma^2\{c_n- y(0)\}-y'(0)\sigma) = \\
%  =\{  (m(m+2R+1)-n(n+2R+1))c_n\},
%\end{multline}
where the formula is valid for  $\lambda_m=m(m+2R+1),$ and $m\in\mathbb{N}$, as in \cite{Kr2}[p. 262].
After simplification in (\ref{LaT1}) it is obtained
\begin{multline}
\label{Lat2N}% \left\{2ny(0)+2y'(0)-2\sum_{k=1}^nkc_{n-k}\right\} =\{  (m(m+2R+1)-n(n+2R+1))c_n\}\\
\left\{2ny(0)+2y'(0)-2\sum_{k=1}^nkc_{n-k}\right\} =\{(m-n)(2R+n+m+1)c_n \}.
\end{multline}
Set $y(x)=\sum_{k=0}^\infty c_k L_k(x)$, given that the Laguerre polynomials satisfy $L_n(0)=1$, $L'_n(0)=-n$ for all $n\geq 0$, \cite[thm.6.3 ]{Be}, then
\begin{eqnarray}
\label{cilag} y(0)=\sum_{k=0}^\infty c_k,\quad y'(0)=-\sum_{k=0}^\infty kc_k
\end{eqnarray}
It is possible to show that solutions of (\ref{La4}) are polynomials of degree $m$. After simplification and taking into account that formulas (\ref{cilag}) are consequently finite,  equation (\ref{Lat2N}) becomes simply
\begin{eqnarray}
\label{goodLag} \{(m-n)(2R+n+m+1)c_n \}=-2\left\{\sum_{k=n+1}^m (k-n)c_k\right\},\quad 0\leq n\leq m.
\end{eqnarray}
Formula (\ref{goodLag}) and conditions $\displaystyle{  y(0)=\sum_{k=0}^m c_k,\quad y'(0)=-\sum_{k=0}^m kc_k}$ provide systems of upper triangular equations,
\begin{align*}
%\label{LagsysN}
%\begin{cases}
c_1+\phantom{000000000}\,c_2+\phantom{00000000000}\cdots\phantom{00000000}+ c_m &=  R\\
(m-1)(2R+m+2)c_1+\phantom{00000000}2(c_2+\phantom{0000}+2c_3+\cdots+(m-1)c_{m}) &= 0\\
(m-2)(2R+m+3)c_2+2(c_3+2c_4+\cdots+(m-2)c_{m})& = 0\\
&\vdots\\
2(R+m)c_{m-1}+2c_m&=0,
%\end{cases}
\end{align*}
 which are easily solved: $c_0=c_1=\cdots =c_{m-1}=-1$, $c_m=R+m$. Therefore for $m=0,1,2,\dots$ it is obtained a family of eigenfunctions $y_m(x),$ which are  the solutions (\ref{solLag}) of the eigenvalue problem (\ref{La4}).

%%%%%%%%%%%%%%%%%%%%%%%%%   conclusions
\section{Conclusions}
%The finite transform method allows for solving differential equations using purely algebraic techniques.
 %This fact opens the possibility of implementing efficient computational algorithms in future work.
%For future work, it is also necessary to expand and generalize the known Laguerre transforms.% as well as to implement as possible Mikusi\'nski's Operational Calculus for Laguerre polynomials. 
The reader may notice that the solutions obtained for the Schr\"{o}\-din\-ger equation in this paper, are slightly different from the approximations already known \cite{G}, \cite{V} in terms of Laguerre polynomials. 
Although the solutions obtained in the present study include a complete series and not just a few terms approximation, the discrepancy depends on the convergence of the Sturm-Liuville expansion of the solutions found. 
The solutions in this article are also different from the solutions in \cite{Al}, however, the reader may note that the polynomials used in \cite{Al} are normalized, which is not the case for the polynomials used here.
In addition, the author in \cite{Al} obtained, by solving the Morse oscillator, a second-order difference equation containing terms in $n^2$.
Terms of order $n^2$ can be obtained also by taking Laguerre transforms of equation (\ref{Moc1}) multiplied by $x$.
But in doing so, the author of this article was unable to reduce the obtained system to a system involving recursive series of only three terms, which is needed essentially in the method in \cite{Al}.
 Perhaps it can be achieved, using equivalent formulas of Laguerre transforms, but this goal was no longer pursued since a solution in terms of Meixner-Polaczek and Laguerre polynomials was obtained, and therefore the goals of this article were achieved.

\section{Appendix}
\label{app}

 In this section are included the proofs of theorem 1 and theorem 2.

%%%%%%%%%%%%    PROF OF THEOREM 2  %%%%%%%%%%%%%%%%
\vspace{5pt}
\noindent {\bf Proof of theorem 1.} Formula (\ref{tnu1})  is shown in \cite{P}[p. 11-16, 11-17]. Formula   (\ref{tnu2}) follows from identity in \cite{Er}[formula (39) p. 192]
\begin{eqnarray*}
L_n^\nu(x)=\sum_{k=0}^n\frac{(\nu-\beta)_k}{k!} L_{n-k}^\beta(x),
\end{eqnarray*}
where $ (r)_k=r(r+1)\cdots (r+k-1)$.

Note that formula (\ref{tnu2}) is obtained from last formula since for $\nu-\beta=1$, it follows $\mathcal{T}^\nu\left[\frac{f(x)}{x}\right]=\left\{\sum_{k=0}^n \frac {(1)_k}{k!}c_{n-k}^{\nu-1}\right\}=\left\{ \sum_{k=0}^nc_{k}^{\nu-1}\right\}$, given that $(1)_k=k!$.

 The proof of formula  (\ref{tnu3}), follows from formula (ii) of theorem 6.11 in~\cite{Be}:
\begin{eqnarray}
\label{nu3} 
xL_n^\nu(x)=(2n+\nu+1)L_n^\nu-(n+\nu)L_{n-1}^\nu(x)-(n+1)L_{n+1}^\nu(x).
\end{eqnarray}
Effectively,
\begin{multline*}
\mathcal{T}^\nu[xf(x)]=\int_0^\infty e^{-x}x^\nu L_n^\nu(x)xf(x) dx =\\
 =\int_0^\infty e^{-x}x^\nu\left[(2n+\nu+1)L_n^\nu-(n+\nu)L_{n-1}^\nu(x)-(n+1)L_{n+1}^\nu(x)\right]f(x)dx=\\
= (2n+\nu+1)c_n^\nu-(n+\nu)c_{n-1}^\nu-(n+1)c_{n+1}^\nu.
\end{multline*}
Formula (\ref{tnu4}) is obtain as follows
\begin{multline*}
 \mathcal{T}^\nu\left[ x\frac{d}{dx}f(x)\right]=\left\{\int_0^\infty e^{-x}x^\nu L_n^\nu(x)x\frac{d}{dx}f(x)dx  \right\}\\
 = \left\{\int_0^\infty e^{-x}x^{\nu+1}\left( L_n^{\nu+1}(x) -L_{n-1}^{\nu+1}(x)\right) \frac{d}{dx}f(x)dx  \right\},
\end{multline*}
where the last integral follows from the identity 
$ L_n^{\nu}(x)=L_n^{\nu+1}(x)-L_{n-1}^{\nu+1}(x),$ \cite{Be}[Theorem 6.11 (i)].
Now, from formula (\ref{tnu1})
\begin{eqnarray*} % ñ
\int_0^\infty e^{-x}x^{\nu+1} L_n^{\nu+1}(x)\frac{d}{dx}f(x)dx=c_n^{\nu+1}-(\nu+1)\sum_{k=0}^nc_{k}^{\nu}+\sum_{k=0}^{n-1}c_k^{\nu+1}
\end{eqnarray*}
\begin{eqnarray*}
\int_0^\infty e^{-x}x^{\nu+1} L_{n-1}^{\nu+1}(x)\frac{d}{dx}f(x)dx=c_{n-1}^{\nu+1}-(\nu+1)\sum_{k=0}^{n-1}c_{k}^{\nu}+\sum_{k=0}^{n-2}c_k^{\nu+1}.
\end{eqnarray*}
In this way, by substracting the last two transforms
\begin{eqnarray*}
\mathcal{T}^\nu\left[ x\frac{d}{dx}f(x)\right]=c_n^{\nu+1}-(\nu+1)c_n^\nu.
\end{eqnarray*}
By using the recurrence relation  $c_n^{\nu+1}=(n+\nu+1)c_n^\nu-(n+1)c_{n+1}^\nu$, \cite{P}[p {\bf 11}-18,11.195 (a)] it is obtained
\begin{eqnarray*}
\mathcal{T}^\nu\left[ x\frac{d}{dx}f(x)\right]=nc_n^\nu-(n+1)c_{n+1}^\nu
\end{eqnarray*}
as claimed.

Finally, formula  (\ref{tnu5}) follows from the solutions of equation (\ref{Laguerre}), equation which is equivalent to
\begin{eqnarray}
\label{Laguerre2} x\frac{d^2}{dx^2}f(x)+(\nu+1-x)\frac{d}{dx}f(x)=-nf(x),
\end{eqnarray}
from (\ref{Laguerre2}) the next transformation follows
\begin{eqnarray}
\label{LaguerreTL} \mathcal{T}^\nu\left[ x\frac{d^2}{dx^2}f(x)+(\nu+1-x)\frac{d}{dx}f(x)\right]=\{-nc_n^\nu\}.
\end{eqnarray}
Since $\mathcal{T}^\nu$ is linear the following relation holds
\begin{eqnarray}
\label{LaguerreTl2} \mathcal{T}^\nu\left[ x\frac{d^2}{dx^2}f(x)\right]=\{-nc_n^\nu\}- \mathcal{T}^\nu\left[ (\nu+1)\frac{d}{dx}f(x)\right]+\mathcal{T}^\nu\left[x\frac{d}{dx}f(x)\right].
\end{eqnarray}
So, formula (\ref{tnu5}) follows directly from (\ref{LaguerreTl2}) and formulas (\ref{tnu3}), and (\ref{tnu4}).\hfill$\square$

\noindent{\bf Remark.} Observe that if the operator $\mathcal{L}^\nu[f(x)]$ is defined as
\begin{eqnarray*}
\mathcal{L}^\nu[f(x)]=e^x x^{-\nu}\frac{d}{dx}\left(x^{\nu+1}e^{-x}\frac{d}{dx}y(x)\right),
\end{eqnarray*}
then formula (\ref{LaguerreTL}) can be obtained from the Sturm-Liouville eigenvalue relation
\begin{eqnarray*}
\mathcal{T}^\nu[\mathcal{L}^\nu[L_n^\nu(x)]]=-n\mathcal{T}^\nu[L_n^\nu(x)],
\end{eqnarray*}
which is probably the most important  Laguerre-type transform  for the purposes of this article.
%%%%%%%%%%%%%%%%%%      END PROOF OF THEOREM 2   %%%%%%%%%%%%%%%%%%%%%%%%%%%%%%

%%%%%%%%%%%%%%%%%%%%%%%%%%   PROOF OF THEOREM 1
\noindent {\bf  Proof of Theorem 2.} Proof by induction on $m$. For $m=0$, $y_0(x)=R=R_0(x)$ is obtained directly from substitution $m=0$ in (\ref{solLag}), taking into account that $L_0(x)=1$ for $x\in [0,\infty)$, as well as direct substitution in (\ref{ksol}). Assume that for $m=1,2,\dots, n$, formula $y_m(x)=R_m(x)$ holds, it is necessary to show that $y_{n+1}(x)=R_{n+1}(x)$.
It follows from formula (\ref{solLag}) that
\begin{multline}
\label{indLa1} y_{n+1}=-L_0(x)-\cdots-L_{n-1}(x)- L_n(x)+(R+n+1)L_{n+1}(x)\\
              =-L_0(x)-\cdots-L_{n-1}(x)- L_n(x)+(R+n+1)L_{n+1}(x)+\\
+(R+n)L_n(x)-(R+n)L_n(x). 
\end{multline}
By substituting the inductive hypothesis in (\ref{indLa1}) it follows that
\begin{eqnarray}
\label{indLa2} y_{n+1}(x)=R_n(x)+(R+n+1)(L_{n+1}(x)-L_n(x))
\end{eqnarray}
But 
\begin{multline*}
L_{n+1}(x)-L_n(x)=\sum_{k=0}^{n+1}\frac{(-1)^k(n+1)!}{(k!)^2(n+1-k)!}x^k-\sum_{k=0}^{n}\frac{(-1)^k n!}{(k!)^2(n-k)!}x^k\\
    =\sum_{k=0}^{n+1}\frac{(-1)^kn!}{k!(k-1)!(n-k+1)!}x^k%+\frac{(-1)^{n+1}}{(n+1)!}x^{n+1}
\end{multline*}
So that,
\begin{multline*}
y_{n+1}(x)=R_n(x)+(R+n+1)\sum_{k=0}^{n+1}\frac{(-1)^kn!}{k!(k-1)!(n-k+1)!}x^k\\
          =R_n(x)+(R+n+1)\sum_{k=0}^{n+1}\frac{(-1)^kn!(k+1)}{(k+1)!(k-1)!(n-k+1)!}x^k\\
  =R_n(x)+\sum_{k=0}^{n+1}\frac{(-1)^k}{(k+1)!}\frac{n![k(R+n+1)+R+n+1]}{(k-1)!(n-k+1)!}x^k\\
 =R_n(x)+\sum_{k=0}^{n+1}\left[\frac{n!(k(R+n+1)+R)}{(k-1)!(n-k+1)!}+\frac{(n+1)!k}{(n-k+1)!k!}\right]x^k\\
 =R_n(x)+\sum_{k=0}^{n+1}\frac{(-1)^k}{(k+1)!}\left[   {n\choose k-1}(k(R+n+1)+R)+{n+1\choose k}k\right]x^k.
\end{multline*}
By using the formula for $R_n(x)$ and the Pascal rule $\displaystyle{{n+1\choose k}={n\choose k}+{ n\choose k-1}}$ it follows that
\begin{multline*}
y_{n+1}(x)=\sum_{k=0}^{n}\frac{(-1)^k}{(k+1)!}  {n\choose k}(k(R+n+1)+R)x^k+\\
+\sum_{k=0}^{n}\frac{(-1)^k}{(k+1)!} {n\choose k-1}(k(R+n+1)+R)x^k+\\
       +\sum_{k=0}^{n+1}\frac{(-1)^k}{(k+1)!}{n+1\choose k}kx^k,
\end{multline*}
consequently,
\begin{multline*}
y_{n+1}(x)   =\sum_{k=0}^{n}\frac{(-1)^k}{(k+1)!}  {n+1\choose k}(k(R+n+1)+R)x^k+ \sum_{k=0}^{n+1}\frac{(-1)^k}{(k+1)!}{n+1\choose k}kx^k\\
=R_{n+1}(x).
\end{multline*}
Therefore $y_m(x)=R_m(x)$ for all $m\in\mathbb{N}$ and $x\in [0,\infty)$.\hfill$\square$

%%%%%%%%%%%%%%%%%%%%%%%%%%%%%%%%%%%   END PROOF TH 1

%%%%%%%%%%%%%%%%%%%%    BIBLIO 
%\bibliographystyle{model1-num-names}

\end{document}